\documentclass[11pt]{article}

\usepackage[T2A]{fontenc}
\usepackage[cp1251]{inputenc}
\usepackage{amsfonts}
\usepackage{eufrak}
\usepackage{amssymb}
\usepackage{amsmath}



\begin{document}

\bigskip

\centerline{\textbf{On a characterisation theorem for probability}}
\centerline{\textbf{distributions on  discrete Abelian groups}}

\bigskip

\centerline{\textbf{G.M. Feldman}}

\bigskip

 \makebox[20mm]{ }\parbox{125mm}{ \small{ Let  $X$ be a countable discrete Abelian group
  containing no elements of order 2,  $\alpha$ be an automorphism of
 $X$,
  $\xi_1$ and  $\xi_2$ be independent random variables with values in the group
       $X$  and distributions
  $\mu_1$ and $\mu_2$. The main result of the article is the following statement.
  The symmetry of the conditional distribution of the linear form
$L_2 = \xi_1 + \alpha\xi_2$ given  $L_1 = \xi_1 +
\xi_2$ implies that   $\mu_j$  are shifts of the Haar distribution
of a finite subgroup of $X$ if and only if
 the automorphism $\alpha$ satisfies the condition ${\rm Ker}(I+\alpha)=\{0\}$.
 This theorem is an analogue for discrete Abelian groups the well-known
 Heyde theorem where Gaussian distribution
on the real line is characterized by the symmetry of the conditional
distribution of one linear form of
independent random variables given another. We also prove
  some generalisations of this theorem.}}

\bigskip


\bigskip


\bigskip

\centerline{\textbf{1. Introduction}}

\bigskip

According to Heyde's theorem the Gaussian distribution
on the real line is characterized by the symmetry of the conditional
distribution of one linear form of
$n$ independent random variables given another  (\cite{He}, see
also \cite[\S\,13.4.1]{KaLiRa}).
Some generalisations of Heyde's theorem, where independent
random variables
$\xi_j$ take values in a locally compact Abelian group and coefficients
of linear forms are topological automorphisms of the group, were studied
in  \cite{Fe2}--\cite{Fe3}, \cite{Fe20bb}--\cite{FeTVP}, \cite{My2}--\cite{MyF}, see also \cite[Chapter VI]{Fe5a} and \cite[Chapter VI]{Fe5}.

Let $X$ be a second countable locally compact Abelian group. We will
consider only such groups,
without mentioning it specifically. In particular,
if  $X$ is a discrete group, it means that
$X$ is countable. Denote by ${\rm Aut}(X)$ the group of topological
automorphisms of a locally compact Abelian group $X$, and by $I$ the identity
automorphism of a group.
The following theorem was proved in \cite{Fe3}, see also \cite[Corollary 17.31]{Fe5}.

{\bf Theorem A.}  {\it  Let $X$ be a discrete Abelian group containing no elements of order $2$.
Let $\alpha$ be an automorphism of the group $X$ satisfying the conditions
 \begin{equation}\label{1a}
I\pm\alpha\in{\rm Aut}(X).
\end{equation}
Let
  $\xi_1$ and  $\xi_2$ be independent random variables with values in the group
       $X$  and distributions
  $\mu_1$ and $\mu_2$. If the conditional distribution of  the linear form  $L_2 = \xi_1 +
\alpha\xi_2$  given $L_1 = \xi_1 + \xi_2$ is symmetric, then
  $\mu_j$   are shifts of the Haar distribution
of a finite subgroup $K$ of the group $X$. Moreover, $\alpha (K)=K$.}

It is well known that Gaussian distributions on discrete Abelian groups   are
degenerated (\cite[Chapter IV]{Pa}). Shifts of the Haar distributions
of  finite subgroups  on such groups   play the role of
Gaussian distributions. For this reason we can consider
Theorem A, as an analogue of Heyde's theorem for discrete Abelian
groups.

The main result of the article, Theorem 1, is the following.
We prove that in fact
Theorem A is valid under a condition which is
significantly weaker  then  $(\ref{1a})$, namely
\begin{equation}\label{1}
{\rm Ker}(I+\alpha)=\{0\}.
\end{equation}
Moreover, condition $(\ref{1})$ can not be relaxed. We prove Theorem 1 in \S 2, and some generalisations of Theorem 1 in \S 3.

Recall some definitions
and agree on notation. For an arbitrary  locally
compact Abelian group  $X$ denote by $Y$ its character
group, and by  $(x,y)$ the value of a character $y \in Y$ at an
element $x \in X$.
 If $K$ is a closed subgroup of $X$, denote by
 $A(Y, K) = \{y \in Y: (x, y) = 1$ \mbox{ for all } $x \in K \}$
its annihilator.
Denote by $b_X$ the subgroup consisting of all compact elements of
the group $X$, and by $c_X$ the connected component of
zero of the group $X$.
If  $\alpha$ is a continuous endomorphism of the group $X$,
then the adjoint  continuous endomorphism $\tilde\alpha$ of the group $Y$ is defined
as follows
   $(x, \tilde\alpha y) = (\alpha x, y)$ for all
 $x \in X$, $y \in Y$. Note that  $\alpha\in {\rm Aut}(X)$
 if and only if $\tilde\alpha\in {\rm Aut}(Y)$.
  Let $K$ be a closed subgroup of the group $X$ and   $\alpha\in {\rm Aut}(X)$. If $\alpha(K)=K$, then the restriction of $\alpha$ to $K$ is a topological automorphism of the group $K$. We denote it by $\alpha_K$.
 Let  $p$ be a prime number. The $p$-torsion subgroup of an Abelian group
$X$ is the set of all elements of $X$ that have order a power of $p$.
Denote by $X_{p}$ the $p$-torsion subgroup of $X$.
Let $n$ be an integer. Denote by $f_n:X \mapsto X$ an endomorphism of
the group $X$, defined by the formula  $f_nx=nx$, $x\in X$.
Put $f_n(X)=X^{(n)}$,
 ${\rm Ker}f_n=X_{(n)}$.
 Denote by $\mathbb{R}$   the group of real numbers.
Let $f(y)$ be a function on the group    $Y$,   and let $h \in
Y$. Denote by   $\Delta_h$   the finite difference operator
$$
\Delta_h f(y) = f(y + h) - f(y), \quad y \in Y.
$$
A function  $f(y)$  on  $Y$ is called a polynomial if
$$
\Delta_h^{n} f(y) = 0
$$
for some natural $n$ and all $y, h\in Y$.

 Denote by ${\rm M}^1(X)$ the
convolution semigroup of probability distributions on the group $X$. Let
${\mu\in {\rm M}^1(X)}$.
Denote by  $$\hat\mu(y) =
\int_{X}(x, y)d \mu(x), \quad y\in Y,$$  the characteristic function of
the distribution  $\mu$, and by $\sigma(\mu)$ the support of $\mu$.
Define the distribution $\bar \mu \in {\rm M}^1(X)$ by the formula
 $\bar \mu(B) = \mu(-B)$ for any Borel subset $B$ in $X$.
Then $\hat{\bar{\mu}}(y)=\overline{\hat\mu(y)}$.
Denote by $m_K$ the Haar distribution of a compact subgroup
 $K$ of the group $X$, and by $E_x$  the degenerate distribution
 concentrated at an element $x\in X$.
We note that the characteristic function of a distribution
$m_K$ is of the form
\begin{equation}\label{11a}
\hat m_K(y)=
\begin{cases}
1, & \text{\ if\ }\   y\in A(Y, K),
\\  0, & \text{\ if\ }\ y\not\in
A(Y, K).
\end{cases}
\end{equation}
Denote by  $\Gamma({\mathbb R}^n )$ the set of Gaussian distributions
on the group ${\mathbb R}^n $.

We will use in the article the standard results on the structure of
locally compact Abelian groups and the duality theory
(see \cite{Hewitt-Ross}).

\bigskip

\centerline{\textbf{2. Proof of the main theorem}}

\bigskip

The main result of the article is the proof of the following statement.

{\bf Theorem 1.}  {\it  Let  $X$ be a  discrete Abelian group
  containing no elements of order $2$. Let  $\alpha$ be an automorphism of
 the group $X$ satisfying condition $(\ref{1})$.
Let
  $\xi_1$ and  $\xi_2$ be independent random variables with values in the group
       $X$  and distributions $\mu_1$ and $\mu_2$.
If the conditional distribution of the linear form
$L_2 = \xi_1 + \alpha\xi_2$ given  $L_1 = \xi_1 +
\xi_2$  is symmetric, then
$\mu_j=m_K*E_{x_j}$, where $K$ is a finite
subgroup of the group $X$,
$x_j\in X$, $j=1, 2$. Moreover, $\alpha (K)=K$.}

To prove Theorem 1 we need some lemmas.

{\bf  Lemma 1} (\cite[Lemma 16.1]{Fe5a}). {\it  Let $X$ be a locally compact Abelian group, $Y$ be
its character group. Let $\alpha$ be a topological automorphism of the group  $X$.
 Let $\xi_1$ and $\xi_2$ be independent random variables with values in the group
       $X$  and distributions
 $\mu_1$ and $\mu_2$.
 The conditional distribution of the linear form
 $L_2 = \xi_1 + \alpha\xi_2$
 given $L_1 = \xi_1 + \xi_2$ is symmetric if and only
 if the characteristic functions
 $\hat\mu_j(y)$ satisfy the equation}
\begin{equation}\label{2a}
\hat\mu_1(u+v )\hat\mu_2(u+\tilde\alpha v )=
\hat\mu_1(u-v )\hat\mu_2(u-\tilde\alpha v), \quad u, v \in Y.
\end{equation}

It is convenient for us to formulate as a lemma the following well-known statement
 (see e.g. \cite[Proposition 2.13]{Fe5a}).

{\bf Lemma 2}. {\it Let $X$ be a locally compact Abelian group,
$Y$ be its character group. Let    $\mu\in{\rm
M}^1(X)$. Then the set  $E=\{y\in Y:\ \hat\mu(y)=1\}$ is a closed
subgroup of the group
   $Y$    and  $\sigma(\mu)\subset A(X,E)$.}

{\bf Lemma 3} (\cite{My1}). {\it Let $X$ be a locally compact Abelian group,  $\alpha$ be a topological automorphism of the group $X$.
Let $\xi_1$ and $\xi_2$ be independent random variables with values in the group
       $X$.
If the   conditional distribution of the linear form
 $L_2 = \xi_1 + \alpha\xi_2$
 given $L_1 = \xi_1 + \xi_2$ is symmetric, then the linear forms
$M_1=(I+\alpha)\xi_1+2\alpha\xi_2$ and
$M_2=2\xi_1+(I+\alpha)\xi_2$ are independent.}

The following lemma is standard and was proved under assumption
 that
$\alpha_j$, $\beta_j$ are topological automorphisms. The proof holds true if
 $\alpha_j$, $\beta_j$ are continuous endomorphisms. We formulate the lemma
 for two independent random variables.

{\bf Lemma 4} (\cite[Lemma 10.1]{Fe5a}). {\it Let $X$ be a locally
compact Abelian group, $Y$ be
its character group. Let $\alpha_j$, $\beta_j$ be continuous
endomorphisms of the group $X$. Let $\xi_1$ and $\xi_2$ be
independent random variables with values in the group $X$  and distributions
 $\mu_1$ and $\mu_2$. The linear forms
 $L_1 = \alpha_1\xi_1 +
\alpha_2\xi_2$ and $L_2 = \beta_1\xi_1  + \beta_2\xi_2$ are independent
if and only
 if the characteristic functions
 $\hat\mu_j(y)$ satisfy the equation}
\begin{equation}\label{3a}
\hat\mu_1(\tilde\alpha_1 u+\tilde\beta_1 v)\hat\mu_2(\tilde\alpha_2 u+\tilde\beta_2 v)=\hat\mu_1(\tilde\alpha_1
u)\hat\mu_2(\tilde\alpha_2
u)\hat\mu_1(\tilde\beta_1 v)\hat\mu_2(\tilde\beta_2 v),
\quad u, v \in
Y.
\end{equation}

{\bf Lemma 5}. {\it Let  $X$ be a  discrete torsion Abelian group,
  containing no elements of order $2$. Let $Y$ be
its character group. Let  $\alpha$ be an automorphism of
the group $X$ satisfying condition $(\ref{1})$.
Let
  $\xi_1$ and  $\xi_2$ be independent random variables with values in the group
       $X$  and distributions $\mu_1$ and $\mu_2$ such that $\hat\mu_j(y) \ge
0$, $j =1, 2$.
If the conditional distribution of the linear form
$L_2 = \xi_1 + \alpha\xi_2$ given  $L_1 = \xi_1 +
\xi_2$  is symmetric, then
  $\hat\mu_1(y)=\hat\mu_2(y)=1$, $y\in B$, where $B$ is an open subgroup of $Y$.}

{\bf Proof.}
By Lemma 3, the symmetry of the conditional distribution of the linear form
 $L_2$ given $L_1$  implies that
the linear forms
$M_1=(I+\alpha)\xi_1+2\alpha\xi_2$ and
$M_2=2\xi_1+(I+\alpha)\xi_2$ are independent.
Then, by  Lemma 4,  the characteristic functions
 $\hat\mu_j(y)$ satisfy equation
$(\ref{3a})$ which takes the form
\begin{equation}\label{4a}
\hat\mu_1((I+\tilde\alpha) u+2 v)\hat\mu_2(2\tilde\alpha  u+(I+\tilde\alpha) v)=\hat\mu_1((I+\tilde\alpha)
u)\hat\mu_2(2\tilde\alpha
u)\hat\mu_1(2 v)\hat\mu_2((I+\tilde\alpha) v), \quad u, v \in
Y.
\end{equation}
Since $X$ is a discrete torsion Abelian group, it follows that $Y$ is
a compact totally disconnected Abelian group.  Then any neighborhood of zero
of the group $Y$ contains an open subgroup.
Hence, we can chose an open subgroup
 $U$   in the group $Y$ such that $\hat\mu_j(y)>0$,
 $y\in U$,
 $j=1, 2$. Put $\psi_j(y)=-\log\hat\mu_j(y)$,    $y\in U$, $j=1, 2$.
 Since $\hat\mu_j(y)\le 1$, it follows that
 \begin{equation}\label{19.04.15.1}
\psi_j(y)\ge 0, \quad y\in U, \ j=1, 2.
\end{equation}
Put
  $V=U\cap\tilde\alpha^{-1}(U)\cap\tilde\alpha^{-2}(U)$. Then $V$ is an open subgroup
  in $Y$. It follows from
$(\ref{4a})$ that the functions $\psi_j(y)$ satisfy the equation
\begin{equation}\label{5a}
\psi_1((I+\tilde\alpha) u+2 v)+\psi_2(2\tilde\alpha  u+(I+\tilde\alpha) v)=P(u)+Q(v),
\quad u, v \in
V,
\end{equation}
where
\begin{equation}\label{9a}
   P(y)=\psi_1((I+\tilde\alpha)
y)+\psi_2(2\tilde\alpha
y), \quad Q(y)=\psi_1(2 y)+\psi_2((I+\tilde\alpha) y), \quad y \in
V.
\end{equation}

To solve equation (\ref{5a}) we use the finite difference method.
Let $h_1$ be an arbitrary element of the group $V$.
Substitute in
 (\ref{5a}) $u+(I+\tilde\alpha) h_1$ for $u$ and $v-2\tilde\alpha  h_1$ for $v$.
 Subtracting equation
  (\ref{5a}) from the resulting equation, we get
  \begin{equation}\label{6a}
    \Delta_{(I-\tilde\alpha)^2 h_1}{\psi_1((I+\tilde\alpha) u+2 v)}
    =\Delta_{(I+\tilde\alpha) h_1} P(u)+\Delta_{-2\tilde\alpha  h_1} Q(v),
\quad u, v\in V.
\end{equation}
Let $h_2$ be an arbitrary element of the group $V$.
Substitute in   (\ref{6a})
  $u+2h_{2}$ for  $u$ and
$v-(I+\tilde\alpha)h_{2}$ for $v$. Subtracting equation
  (\ref{6a}) from the resulting equation, we obtain
 \begin{equation}\label{7a}
     \Delta_{2 h_2}\Delta_{(I+\tilde\alpha) h_1} P(u)+\Delta_{-(I+\tilde\alpha) h_2}\Delta_{-2\tilde\alpha  h_1} Q(v)=0,
\quad u, v\in V.
\end{equation}
Let $h$ be an arbitrary element of the group $V$.
Substitute in
(\ref{7a})   $u+h$ for $u$. Subtracting equation
  (\ref{7a}) from the resulting equation, we have
 \begin{equation}\label{8a}
   \Delta_{h}\Delta_{2 h_2}\Delta_{(I+\tilde\alpha) h_1} P(u)=0,
\quad u\in V.
\end{equation}

 Put $W_P=V\cap V^{(2)}\cap(I+\tilde\alpha)(V)$. Since
$h, h_1, h_2$ in (\ref{8a})
 are arbitrary elements of the group
  $V$, it follows from
 (\ref{8a}) that
   the function $P(y)$ satisfies the equation
\begin{equation}\label{30_07_1}
\Delta_h^{3} P(y) = 0, \quad y, h  \in W_P.
\end{equation}

Put $W_Q=V\cap (I+\tilde\alpha)(V)\cap (\tilde\alpha(V))^{(2)}$.
Reasoning similarly, we obtain from
(\ref{7a}) that the function $Q(y)$ satisfies the equation
\begin{equation}\label{30_07_2}
\Delta_h^{3} Q(y) = 0, \quad y, h  \in W_Q.
\end{equation}  Set $W=W_P\cap W_Q$.

Since the automorphism $\alpha$
  satisfies condition  $(\ref{1})$,
 the group $(I+\tilde\alpha)(Y)$ is dense in $Y$, and taking into account that
  $Y$ is a compact group, we have
  \begin{equation}\label{2b}
  (I+\tilde\alpha)(Y)=Y.
\end{equation}
It follows from this that the continuous endomorphism
     $I+\tilde\alpha$
 is open, and hence $(I+\tilde\alpha)(V)$ is an open subgroup
 in $Y$. Since
  $X$ is a torsion group containing no elements of
  order 2, we have $f_2\in {\rm Aut}(X)$.
It follows from this that $f_2\in {\rm Aut}(Y)$, and hence
 $V^{(2)}$ and $(\tilde\alpha(V))^{(2)}$ are open subgroups in $Y$.
 This implies that $W$ is an open subgroup
 in $Y$. Hence, $W$ is a closed subgroup
 in $Y$, so that $W$ is a compact subgroup. It follows from
 (\ref{30_07_1}) and (\ref{30_07_2}) that $P(y)$ and $Q(y)$ are continuous polynomials on the group $W$.
 It is well known, see e.g. \cite[Proposition 5.7]{Fe5a}, that a
continuous polynomial on a compact Abelian group is a constant. Since
 $P(0)=Q(0)=0$, we have
 \begin{equation}\label{10a}
   P(y)=Q(y)=0, \quad y\in W.
\end{equation}
Put $B=(I+\tilde\alpha)(W)$. It is obvious that
$B\subset U$.
Since $I+\tilde\alpha$ is a continuous open endomorphism,
$B$ is an open subgroup in $Y$.
It follows from $(\ref{19.04.15.1})$, $(\ref{9a})$
and $(\ref{10a})$ that $\psi_1(y)=\psi_2(y)=0$, $y\in B$.
Hence, $\hat\mu_1(y)=\hat\mu_1(y)=1$,  $y\in B$. Lemma 5 is proved.

It should be noted that if an automorphism
 $\alpha$ satisfies conditions
 (\ref{1a}), then Lemma 5 is a particular case
 of a general statement which was proved
 in \cite[Lemma 6]{Fe3}, see also \cite[Lemma 17.24]{Fe5a}.

It is convenient for us to formulate as a lemma the following easily verified statement.

{\bf Lemma 6.} {\it  Let $X$ be a locally
compact Abelian group, $Y$ be
its character group.  Let $G$ be a compact subgroup
of a group $X$ and $\beta$ be a continuous endomorphism of the
group $X$. Then the following statements are equivalent:

$(i)$ $\beta(G) \supset G;$

$(ii)$ if $\tilde\beta y \in A(Y, G)$, then $y \in A(Y, G)$.}

{\bf Lemma 7.} {\it Let $X$ be a locally
compact Abelian group containing no elements of order $2$.
Let $\alpha$ be a topological automorphism of the group  $X$.
Let $\xi_1$ and $\xi_2$ be independent random variables with
values in a group  $X$ and  distributions $\mu_1=m_{K_1}$ and $
\mu_2=m_{K_2}$, where $K_1$ and $K_2$ are finite subgroups of $X$.
 If the
conditional distribution of the linear form  $L_2 =
 \xi_1 + \alpha\xi_2$ given $L_1= \xi_1+\xi_2$
is symmetric, then $K_1=K_2=K$ and $\alpha(K)=K$.}

 {\bf Proof}. Denote by $Y$ the character group of the group $X$.
 Put $f(y)=\hat m_{K_1}(y)$, $g(y)=\hat m_{K_2}(y)$, $E_j=A(Y, K_j)$,
 $j=1, 2$.
 It follows from (\ref{11a}) that
 \begin{equation}\label{7}
 f(y) = \begin{cases}1, & \text{\ if\ }y\in E_1,\\ 0, &   \text{\ if\ }y\notin E_1,
\\
\end{cases} \quad\quad\quad
g(y) = \begin{cases}1, & \text{\ if\ }y\in E_2,\\ 0, &\text{\ if\ }   y\notin E_2.
\\
\end{cases}
\end{equation}
By Lemma 1, the characteristic functions $f(y)$ and $g(y)$ satisfy
equation
(\ref{2a}) which takes the form
 \begin{equation}\label{2}
f(u+v)g(u+\tilde\alpha v)=f(u-v)g(u-\tilde\alpha v), \quad u, v\in Y.
\end{equation}
 Put in $(\ref{2})$ $u=v=y$. We get
 \begin{equation}\label{3}
f(2y)g((I+\tilde\alpha) y)=g((I-\tilde\alpha) y), \quad  y\in Y.
\end{equation}

Assume that
\begin{equation}\label{4}
(I-\tilde\alpha) y\in E_2.
\end{equation}
Then, it follows from $(\ref{7})$  and $(\ref{3})$ that
\begin{equation}\label{5}
2y\in E_1
\end{equation}
and
\begin{equation}\label{6}
(I+\tilde\alpha) y\in E_2.
\end{equation}
Since group $X$ contains no elements of order $2$ and $K_2$ is a finite group, we have
\begin{equation}\label{14}
(K_2)^{(2)}=K_2.
\end{equation}
It follows from $(\ref{4})$  and $(\ref{6})$ that
\begin{equation}\label{11.08.15}
2y\in E_2.
\end{equation}
Taking into account $(\ref{14})$ and applying Lemma 6 to
$\beta=f_2$, $G=K_2$, we  get from $(\ref{11.08.15})$ that $y\in E_2$. So, we proved that  $(\ref{4})$ implies that $v\in E_2$. Applying Lemma 6 again to $\beta=I-\alpha$, $G=K_2$, we obtain that
$(I-\alpha)(K_2)\supset K_2$.
Since $K_2$ is a finite group, it follows from this that
\begin{equation}\label{17}
(I-\alpha)(K_2)=K_2,
\end{equation}
and hence,
\begin{equation}\label{9}
\alpha(K_2)=K_2.
\end{equation}

Assume that $y\in E_2$. It is obvious that $(\ref{17})$ implies the inclusion $(I-\tilde\alpha)(E_2)\subset  E_2$. Hence,
 $(\ref{4})$ holds, and then $(\ref{5})$ follows from $(\ref{7})$ and $(\ref{3})$.
 Since   group $X$ contains no elements of order $2$ and $K_1$ is a finite group, we have
$(K_1)^{(2)}=K_1$.
Taking into account $(\ref{5})$ and applying Lemma 6
to  $\beta=f_2$, $G=K_1$, it follows from $(K_1)^{(2)}=K_1$
that $y\in E_1$. So,
we proved that
\begin{equation}\label{12}
E_2\subset E_1.
\end{equation}

Put in   $(\ref{2})$ $u=\alpha y$, $v=y$. We get
 \begin{equation}\label{10}
f((I+\tilde\alpha)y)g(2\tilde\alpha y)=f((I-\tilde\alpha) y), \quad  y\in Y.
\end{equation}
Arguing similarly as above we get that $(\ref{10})$ implies the inclusion
 \begin{equation}\label{13_08_17_1}
{(I-\tilde\alpha)(E_1)\subset E_1}.
\end{equation}

Let $y\in E_1$. It follows from $(\ref{13_08_17_1})$  that  $(I-\tilde\alpha) y\in E_1$, and we get from $(\ref{7})$ and $(\ref{10})$ that
 \begin{equation}\label{15}
2\tilde\alpha y\in E_2.
\end{equation}
It is easy to see  $(\ref{14})$  and $(\ref{9})$ imply
$2\alpha(K_2)=K_2$.
Taking into account $(\ref{15})$    and applying Lemma 6
to
$\beta=2\alpha$, $G=K_2$, it follows from $2\alpha(K_2)=K_2$
that $y\in E_2$.

So, we proved that
$E_1 \subset E_2$.
Considering  $(\ref{12})$, this implies that $E_1=E_2$, and hence,
$K_1=K_2=K$. The equality $\alpha(K)=K$ follows from $(\ref{9})$. Lemma 7 is proved.

Note that if $f_2\in {\rm Aut}(X)$ and a topological  automorphism
 $\alpha$ satisfies conditions
 (\ref{1a}), then   Lemma 7 was proved
 in \cite[Lemma 7]{Fe3}, see also \cite[Lemma 17.25]{Fe5a}.

  {\bf Remark 1}. Generally speaking, Lemma 7 fails if a group $X$ contains an element  of order $2$.
 Denote by $\mathbb{Z}(2)=\{0, 1\}$ the group of residues modulo 2. Let $X=(\mathbb{Z}(2))^2$. Denote by $x=(x_1, x_2)$, where  $x_j\in \mathbb{Z}(2)$, elements of the group $X$.
Let $\xi_1$ and  $\xi_2$ be independent   random variables with values in the group $X$. It follows from Lemma 1 that for any automorphism $\alpha$
the conditional distribution of the linear form
$L_2 = \xi_1 + \alpha\xi_2$ given  $L_1 = \xi_1 +
\xi_2$  is symmetric.
Put $K_1=\{(x_1, 0)\}$, $K_2=\{(0, x_2)\}$, $x_j\in \mathbb{Z}(2)$. Let $\xi_1$ and  $\xi_2$ be independent   random variables with values in the group $X$ and distributions $m_{K_1}$ and $m_{K_2}$. Then the conditional distribution of the linear form
$L_2 = \xi_1 + \alpha\xi_2$ given  $L_1 = \xi_1 +
\xi_2$  is symmetric, whereas $K_1\ne K_2$.

 {\bf Remark 2}. Generally speaking, Lemma 7 fails if $K_1$ and $K_2$ are
compact but not finite subgroups of the group $X$.
Here is an example.
(compare  with  \cite[Remark 13.16]{Fe5a}). Let $G$ be an arbitrary compact Abelian
group. Consider the direct product  $$X=\mathop{\mbox{\rm\bf
P}}\limits_{j\in \mathbb{Z}} G_j,$$
where all $G_j = G$.
  Put
$H=G^*$. Then the group  $Y$ is topologically isomorphic
to a weak direct product of groups
$H_j$, where $H_j=H,$
  $$Y\cong\mathop{\mbox{\rm\bf
P}^*}\limits_{j\in \mathbb{Z}} H_j.$$
 Let $\alpha \in
{\rm Aut}(X)$ be an automorphism of the form
$$\alpha(g_j)=
(g_{j-2}), \quad\quad
(g_j) \in X.$$
 Then
$$\tilde\alpha(h_j)=(h_{j+2}),
\quad\quad (h_j)\in Y.$$
It is obvious that $${\rm Ker}(I-\alpha)=\{(g_j)\in X: g_{2k}=g_2, \ g_{2k-1}=g_1, \
g_1, g_2\in G, \ k\in \mathbb{Z}\}.$$ Since $(I-\tilde\alpha)(Y)=A(Y, {\rm Ker}(I-\alpha))$, we have
$$(I-\tilde\alpha)(Y)=\{(h_j)\in Y: \sum\limits_{k\in \mathbb{Z}}h_{2k}=0, \ \sum\limits_{k\in \mathbb{Z}}h_{2k-1}=0, \ k\in \mathbb{Z}\}.$$
 Consider the subgroups
$$K_1=\mathop{\mbox{\rm\bf P}}\limits_{j\in \mathbb{Z}, j\neq 1} G_j, \qquad
K_2=\mathop{\mbox{\rm\bf P}}\limits_{j\in \mathbb{Z}, j\neq 2} G_j.$$
 Obviously, $K_1 \ne
K_2$. Put $L_j =A(Y,K_j)$, $j=1, 2$. Denote by $L$ the subgroup of   $Y$, generated by
 $L_1$ and $L_1$. It is easy to see that $L\cap(I-\tilde\alpha)(Y)=\{0\}$ and $L_1\cap L_2=\{0\}$.

 We will check that the
characteristic functions $f(y)=\hat m_{K_1}(y)$ and $g(y)=\hat m_{K_2}(y)$
satisfy equation (\ref{2}). Let the left-hand side of (\ref{2}) is equal to 1. Then
$u+v\in L_1$, $u+\tilde\alpha v\in L_2$. It follows from this that $(I-\tilde\alpha)v\in L$.
Hence, $(I-\tilde\alpha)v=0$. Since, obviously, ${\rm Ker}(I-\tilde\alpha)=\{0\}$,
we have $v=0$. This implies that $u\in  L_1\cap L_2$, and hence, $u=0$.
So, the right-hand side of (\ref{2}) is also equal to  1.
Similarly we verify that if the right-hand side of (\ref{2}) is also equal to  1,
then the left-hand side of (\ref{2}) is also equal to 1.  It means that the
characteristic functions $f(y)$ and $g(y)$ satisfy equation (\ref{2}).

Let $\xi_1$ and
$\xi_2$ be independent random variables with values in the group $X$ and
distributions  $m_{K_1}$ and $m_{K_2}$. By Lemma 1,
the
conditional distribution of the linear form  $L_2 =
 \xi_1 + \alpha\xi_2$ given $L_1= \xi_1+\xi_2$
is symmetric.

 {\bf Lemma 8} (\cite{Fe2}, see also \cite[Corollary 17.2 and Remark 17.5]{Fe5a}).
 {\it  Let  $X$ be a  finite Abelian group
  containing no elements of order $2$. Let  $\alpha$ be an automorphism of
 $X$ satisfying condition $(\ref{1})$.
Let
  $\xi_1$ and  $\xi_2$ be independent random variables with values in the group
       $X$  and distributions $\mu_1$ and $\mu_2$.
If the conditional distribution of the linear form
$L_2 = \xi_1 + \alpha\xi_2$ given  $L_1 = \xi_1 +
\xi_2$  is symmetric, then
$\mu_j=m_K*E_{x_j}$, where $K$ is a
subgroup of the group $X$,
$x_j\in X$, $j=1, 2$. Moreover, $\alpha (K)=K$.}

{\bf Lemma 9} (\cite{FeTVP}). {\it  Let
 $Y$ be a connected compact Abelian group, $X$ be its character group.
 Let $\tilde\alpha$ be a topological automorphism of the group
 $Y$ satisfying condition$ (\ref{2b})$.
Let $\mu_1$ and $\mu_2$ be distributions on the group $X$ such that their
characteristic functions $\hat\mu_j(y)$ satisfy equation  $(\ref{2a})$.
  Then $\hat\mu_j(y)=(x_j, y)$, where  $x_j\in X$, $j=1, 2$.}

  {\bf Lemma 10.}  {\it  Let  $X$ be a  discrete Abelian group, $\alpha$ be an automorphism of
 the group $X$ satisfying condition
 \begin{equation}\label{02.11.1}
{\rm Ker}(I+\alpha)\subset b_X.
\end{equation}
Let
  $\xi_1$ and  $\xi_2$ be independent random variables with values in the group
       $X$  and distributions $\mu_1$ and $\mu_2$.
If the conditional distribution of the linear form
$L_2 = \xi_1 + \alpha\xi_2$ given  $L_1 = \xi_1 +
\xi_2$  is symmetric, then there exist elements
  $x'_1, x'_2\in X$ such that the distributions
$\mu'_j$ of the random variables  $\xi_j'=\xi_j-x'_j$ are supported in the subgroup $b_X$,
and the conditional distribution of the linear form
$L'_2 = \xi'_1 + \alpha\xi'_2$ given $L'_1 = \xi'_1 + \xi'_2$ is symmetric.}

{\bf Proof.} Denote by $Y$ the character group of the group $X$.
Since $X$ is a discrete group, $b_X$ is the subgroup consisting of
all elements of finite order of the group $X$.   Consider the factor-group
$X/b_X$ and denote by $[x]$ its elements.
Note that the character group of the group
$X/b_X$ is topologically isomorphic to the annihilator $A(Y, b_X)$, and
 $A(Y, b_X)=c_Y$. Note that $c_Y$  is a compact group.
It follows from
\begin{equation}\label{4c}
\alpha(b_X)=b_X
\end{equation}
 that $\alpha$
induces an automorphism $\hat\alpha$ of the factor-group $X/b_X$ by the formula
 $\hat\alpha [x]=[\alpha x]$. We note that $\tilde\alpha(c_Y)=c_Y$, and a topological automorphism
 of the group $c_Y$ which is adjoint to $\hat\alpha$ is a restriction of
 $\tilde\alpha$ to $c_Y$.  Verify that
\begin{equation}\label{1b}
(I+\tilde\alpha_{c_Y})(c_Y)=c_Y.
\end{equation}
Taking into account that for a discrete Abelian group $X$ conditions (\ref{1}) and (\ref{2b})
are equivalent, it suffices to verify that
\begin{equation}\label{3b}
{\rm Ker}(I+\hat\alpha)=\{0\}.
\end{equation}

Let $x_0\in X$ and $[x_0]\in {\rm Ker}(I+\hat\alpha)$. Then $[(I+\alpha)x_0]=0$, and this implies that $(I+\alpha)x_0\in b_X$, i.e. $n(I+\alpha)x_0=0$ for some integer $n$. Hence,  $(I+\alpha)nx_0=0$. Taking into account  (\ref{02.11.1}) this implies   that $nx_0\in b_X$. Thus,  $x_0\in b_X$.
 So, (\ref{3b}) holds, thus (\ref{1b}) is proved.

Consider the restriction of equation   (\ref{2a})
to $c_Y$. Taking into account that by Lemma 1, the characteristic functions
$\hat\mu_j(y)$ satisfy equation $(\ref{2a})$, and applying Lemma  9 to the group
 $c_Y$, we get that the restrictions of the characteristic functions
$\hat\mu_j(y)$ to  the subgroup $c_Y$ are characters of the subgroup $c_Y$.
Extending these characters to characters of the group $Y$, we obtain that
there exist elements $x_j\in X$, $j=1, 2$, such that
\begin{equation}\label{18.04.1}
\hat\mu_j(y)=(x_j, y), \quad y\in c_Y, \ j=1, 2
\end{equation}
holds. Substituting  (\ref{18.04.1}) into (\ref{2a}) and taking into account that
$A(X, c_Y)=b_X$, we get
\begin{equation}\label{18.04.2}
2(x_1+\alpha x_2)\in b_X.
\end{equation}
Since $b_X$ consists of all elements of finite order of
the group $X$, it follows from
 (\ref{18.04.2}) that
\begin{equation}\label{18.04.3}
x_1+\alpha x_2\in b_X.
\end{equation}

Put $x_1'=-\alpha x_1$, $x_2'=x_2$.
It is easy to see that (\ref{18.04.1}) and (\ref{18.04.3}) imply that
\begin{equation}\label{18.04.4}
\hat\mu_j(y)=(x'_j, y), \quad y\in c_Y, \ j=1, 2.
\end{equation}
Since $x'_1+\alpha x'_2=0$, the characteristic functions  $f_j(y)=(-x'_j, y)$,
$j=1, 2,$ on the group $Y$ satisfy equation (\ref{2a}). Put $\xi_j'=\xi_j-x'_j$, and
denote by $\mu'_j$ the distribution of the random varisble $\xi_j'$. It follows from
\begin{equation}\label{18.04.5}
\mu'_j=\mu_j*E_{-x'_j}, \quad j=1, 2,
\end{equation}
   that the characteristic functions
$\hat\mu'_j(y)$ also satisfy equation (\ref{2a}). By Lemma 1,
the conditional distribution of the linear form
$L'_2 = \xi'_1 + \alpha\xi'_2$ given
 $L'_1 = \xi'_1 + \xi'_2$   is symmetric. Moreover,  (\ref{18.04.4}) and
(\ref{18.04.5}) imply that $\hat\mu'_j(y)=1$, $y\in c_Y$. Then, by Lemma 2,
 $\sigma(\mu'_j)\subset A(X, c_Y)=b_X$. Lemma 10 is proved.

Note  that if an   automorphism
 $\alpha$ satisfies conditions
 (\ref{1a}), then   Lemma 10 is a particular case of a general
 statement which was proved
 in \cite[Corollary 1]{Fe3}, see also
 \cite[Lemma 17.22]{Fe5a}.

{\bf Remark 3}. Generally speaking, Lemma 10 fails if we omit condition (\ref{02.11.1}).
 Put $G={\rm Ker}(I+\alpha)$. Let $x_0\in G$, $x_0\notin b_X$ and let
  $\mu$ be a distribution on $X$ such that $\sigma(\mu)=\{0, x_0\}$.
 Let $\xi_1$ and $\xi_2$ be
 independent identically distributed random variables with values in the group
  $G$ and distribution $\mu$.
It is obvious that  $\alpha x=-x$ for all
 $x\in G$. Hence, applying  Lemma 1 to the group $G$ we get that
 the conditional distribution of the linear form  $\tilde L_2 = \xi_1 - \xi_2$ given
 $\tilde L_1 = \xi_1 + \xi_2$ is symmetric.
Hence, if we consider
  the random variables   $\xi_1$ and $\xi_2$, as
  random variables with values in the group $X$, then
the conditional distribution of the linear form    $L_2 = \xi_1 +
\alpha\xi_2$ given $L_1 = \xi_1 + \xi_2$ is also symmetric.
 Since $\sigma(\mu*E_x)=\{x, x_0+x\}$ for any $x\in X$, it is obvious that
 does  not exist  an element $x\in X$ such that the distribution $\mu'$ of
 the random variables $\xi_j'=\xi_j+x$ is supported in the subgroup $b_X$.

{\bf Proof of Theorem 1.} We will follow the scheme of the proof
of Theorem 17.26 in \cite{Fe5a}. Denote by $Y$ the character group of the group $X$. Since $\alpha(b_X)=b_X$, we can apply Lemma 10 and
 assume that $X$ is a torsion group. By Lemma 1,
the symmetry of the conditional distribution of the linear form
$L_2$ given $L_1$ implies that the characteristic functions
 $\hat\mu_j(y)$ satisfy equation
(\ref{2a}).
Put
$\nu_j = \mu_j
* \bar \mu_j$. Тогда
 $\hat \nu_j(y) = |\hat \mu_j(y)|^2 \ge 0,$   $y \in Y$.
Obviously,   the characteristic functions $\hat \nu_j(y)$
also satisfy equation (\ref{2a}). If we prove that
 $\nu_1=\nu_2=m_K$, where $K$ is a finite subgroup of the group
  $X$, this easily implies that
$\mu_j=m_K*E_{x_j}$, where $x_j\in X$, $j=1, 2$.
Therefore we can solve  equation (\ref{2a}) supposing that
 $\hat\mu_j(y)\ge 0.$ So, we will assume that.

Put $f(y) = \hat\mu_1(y)$, $g(y) =
\hat\mu_2(y)$. Then equation (\ref{2a}) takes the form (\ref{2}).
We will prove that in this case $f(y) = g(y)=\hat
m_K(y)$,
 where $K$ is a finite subgroup of the group $X$.

Set $E_f= \{y \in Y:
f(y) = 1\}$, $E_g= \{y \in Y: g(y) = 1\}$. By Lemma 2,  $E_f$ and $E_g$ are closed subgroups of
$Y$. Put $F=A(X, E_f)$, $G=A(X, E_g)$.  Then, by Lemma 2,
$\sigma(\mu_1) \subset   F$, $\sigma(\mu_2)
\subset   G$. By Lemma 5, there exists an open subgroup
 $B$ such that $B \subset E_f\cap E_g$. Put $S
= A(X, B)$. Then $F$ and $G$ are subgroups in $S$. Since $B$ is an open subgroup,
 $S$ is a compact subgroup, and taking into account that $X$ is a discrete group,
 $S$ is a finite subgroup. Hence,
 $F$ and $G$ are also finite subgroups.

It follows from   (\ref{2}) that
\begin{equation}\label{12a}
f^n(u+v) g^n(u+\tilde\alpha v)= f^n(u-v) g^n(u-\tilde\alpha v), \quad u,v
\in Y,
\end{equation}
holds for any natural $n$.
It is obvious that there exist limits
\begin{equation}\label{5e}
\bar f(y)=\lim_{n \rightarrow \infty} f^n(y)=
\begin{cases}
1, & \text{\ if\ }\ y\in E_f,
\\  0, & \text{\ if\ }\ y\not\in
E_f,
\end{cases}\quad\quad
\bar g(y)=\lim_{n \rightarrow \infty} g^n(y)=
\begin{cases}
1, & \text{\ if\ }\ y\in E_g,
\\  0, & \text{\ if\ }\ y\not\in E_g.
\end{cases}
\end{equation}
Since   $E_f=A(Y, F)$, $E_g=A(Y, G)$, it follows from
 (\ref{11a})  that
$$
\hat m_F(y)=
\begin{cases}
1, & \text{\ if\ }\ y\in E_f,
\\  0, & \text{\ if\ }\ y\not\in
E_f,
\end{cases}\qquad\hat m_G(y)=
\begin{cases}
1, & \text{\ if\ }\ y\in E_g,
\\  0, & \text{\ if\ }\ y\not\in E_g.
\end{cases}
$$
Hence,
$$\hat m_F(y) = \bar f(y),
\quad \hat m_G(y) = \bar g(y).
$$
We note that since $X$ is a discrete torsion Abelian group containing no elements
of order $2$, we have $f_2\in {\rm Aut}(X)$. Let $\zeta_1$    and
$\zeta_2$  be independent random variables with values in
the group $X$  and distributions $m_F$ and
$m_G$. It follows from (\ref{12a}) and (\ref{5e}) that the characteristic functions
 $\bar f(y)$ and $\bar g(y)$ also satisfy equation
 (\ref{2}).  By Lemma 1, this implies that
the conditional distribution of the linear form
$L_2=\zeta_1+\alpha\zeta_2$ given
$L_1=\zeta_1+\zeta_2$ is symmetric. Applying Lemma 7  we obtain that
 $F=G$ and $\alpha(F)=F$.

We now return to the random variables $\xi_1$  and $\xi_2$ and the linear forms
  $L_1 = \xi_1 + \xi_2$ and $L_2 = \xi_1 + \alpha\xi_2$.
  Since  $\sigma(\mu_j)\subset F$, the random variables $\xi_j$ take values
  in the finite group $F$. Taking into account that $\alpha(F)=F$, we can apply
  Lemma 8 to the group $F$.  Since $\hat\mu_j(y)\ge 0$, by Lemma 8,
 $\mu_1=\mu_2=m_K$, where $K$ is a finite subgroup of the group
 $X$, and $\alpha (K)=K$.   Theorem 1 is proved.

{\bf Remark 4}.   Let $X$ be a discrete Abelian group containing no elements
of order $2$, and $\alpha$ be an automorphism of
 the group $X$.
 Assume that $G={\rm Ker}(I+\alpha)\ne\{0\}$. Let $\xi_1$ and $\xi_2$ be
 independent identically distributed random variables with values in the group
  $G$ and distribution $\mu$.
It is obvious that  $\alpha x=-x$ for all
 $x\in G$. Hence, applying  Lemma 1 to the group $G$ we get that
 the conditional distribution of the linear form  $\tilde L_2 = \xi_1 - \xi_2$ given
 $\tilde L_1 = \xi_1 + \xi_2$ is symmetric.
It follows from this that if we consider
  the random variables   $\xi_1$ and $\xi_2$, as
  random variables with values in the group $X$, then
the conditional distribution of the linear form  $L_2 = \xi_1 +
\alpha\xi_2$ given $L_1 = \xi_1 + \xi_2$ is also symmetric. Taking into account that
$\mu$ is an arbitrary distribution, we see that Theorem 1 fails if condition
(\ref{1}) is not fulfilled.

We complement Theorem 1 by the following statement.

{\bf Proposition 1.}  {\it  Let  $X$ be a  locally compact Abelian group
  containing no elements of order $2$, $K$ be a finite subgroup of $X$,
  $\alpha$ be a topological automorphism of
 $X$.
Let
  $\xi_1$ and  $\xi_2$ be independent identically distributed
  random variables with values in the group
       $X$  and distribution $m_K$. Then the following statements are equivalent:

$(i)$ the conditional distribution of the linear form
$L_2 = \xi_1 + \alpha\xi_2$ given  $L_1 = \xi_1 +
\xi_2$  is symmetric;

$(ii)$ $(I-\alpha )(K)=K$.}

{\bf Proof}. Denote by $Y$ the character group of the group $X$.
Put $L=A(Y, K)$, $f(y)=\hat m_K(y)$.

$(i)\Rightarrow (ii)$. By Lemma 7, it follows from $(i)$ that $\alpha(K)=K$. It means that we can assume that   $X=K$, and the characteristic function $f(y)$ is of the form
\begin{equation}\label{13_08_17_2}
f(y)=
\begin{cases}
1, & \text{\ if\ }\   y=0,
\\  0, & \text{\ if\ }\ y\ne 0.
\end{cases}
\end{equation}
By Lemma 1, the characteristic function
  $f(y)$ satisfies equation
\begin{equation}\label{03.11.1}
f(u+v) f(u+\tilde\alpha v)= f(u-v) f(u-\tilde\alpha v), \quad u,v
\in Y.
\end{equation}
Put in (\ref{03.11.1}) $u=v=y$. We get
\begin{equation}\label{03.11.2}
f(2y) f((I+\tilde\alpha)y)= f((I-\tilde\alpha)y), \quad y
\in Y.
\end{equation}
Let $y\in {\rm Ker} (I-\tilde\alpha)$. Then it follows from  $(\ref{03.11.2})$ that
  $2y=0$. Since $X$ is a finite group and contains no elements of order 2, the group $Y$ also contains no elements of order 2, and hence  $y=0$, i.e. ${\rm Ker} (I-\tilde\alpha)=\{0\}$. This implies   $(ii)$.

$(ii)\Rightarrow (i)$. The characteristic function $f(y)$ is of the form
 (\ref{11a}).  We shall verify that
  $f(y)$ satisfies equation (\ref{03.11.1}). Assume that for some
$u, v \in Y$ the left hand-side of equation (\ref{03.11.1}) is equal to 1.
Then $u+v, u+\tilde\alpha v\in L$.
This implies that $(I-\tilde\alpha )v\in L$. Taking into account
 $(ii)$, by Lemma 6, applying to $\beta=I-\alpha$,
we get $v\in L$. Hence, $u\in L$, and then $\tilde\alpha  v\in L$. It follows from this that
$u-v, u-\tilde\alpha v\in L$. We get that the right hand-side of equation
(\ref{03.11.1}) is also equal to 1.
Reasoning similarly we check that if the right hand-side of equation
(\ref{03.11.1}) is   equal to 1, than the left hand-side of equation (\ref{03.11.1})
is equal to 1. Thus  the characteristic function
  $f(y)$ satisfies equation (\ref{03.11.1}). By Lemma 1,  $(i)$ holds.
   Proposition 1 is proved.

  Note that the proof of the statement $(ii)\Rightarrow (i)$ is based on
   Lemmas 1 and 6 only. Therefore, the statement $(ii)\Rightarrow (i)$
   also holds in the case, when  $K$ is a compact subgroup of $X$.

{\bf Remark 5}. Generally speaking, Proposition 1 fails if a group $X$ contains elements of order 2.
Let a group $X$ and subgroups $K_j$ be as in Remark 1.
Put $\alpha (x_1, x_2)=(x_2, x_1+x_2)$. It is obvious that $\alpha\in {\rm Aut}(X)$ and $\alpha$ satisfies condition $(\ref{1})$. We have $(I-\alpha)K_2=K_1$. Let
  $\xi_1$ and  $\xi_2$ be independent identically distributed
  random variables with values in the group
       $X$  and distribution $m_{K_2}$.
        Then condition $(i)$ holds, whereas  $(ii)$ is not.
\bigskip

\centerline{\textbf{3. Generalizations of Theorem 1}}

\bigskip

First we prove the following generalization of Theorem 1.

{\bf Theorem 2.} {\it  Let $X={\mathbb R}^n\times G$, where $n \ge 0,$ and $G$ is
a discrete Abelian group
  containing no elements of order $2$. Let  $\alpha$ be a topological automorphism of
 the group $X$ satisfying condition $(\ref{1})$.
 Let
  $\xi_1$ and  $\xi_2$ be independent random variables with values in the group
       $X$  and distributions $\mu_1$ and $\mu_2$.
If the conditional distribution of the linear form
$L_2 = \xi_1 + \alpha\xi_2$ given  $L_1 = \xi_1 +
\xi_2$  is symmetric, then
$\mu_j=\gamma_j*m_K*E_{g_j}$, where $\gamma_j \in \Gamma({\mathbb R}^n)$,
 $K$   is a finite
subgroup of the group  $G$,
$g_j\in G$, $j=1, 2$. Moreover, $\alpha (K)=K$.}

To prove Theorem 2 we need the following lemma.

{\bf Lemma 11} (\cite{Fe14}, see also \cite[Theorem 13.3]{Fe5a}). \textit{
 Let $X=\mathbb{R}^n\times G$, where $n\ge 0$, and $G$
is a finite Abelian group. Let   $\xi_1$ and $\xi_2$ be
independent random variables with values in the group  $X$ and
distributions $\mu_1$ and $\mu_2$. Let $\alpha_j, \beta_j$, $j=1,
2,$ be topological automorphisms of the group $X$. If the linear forms
$L_1 = \alpha_1\xi_1 + \alpha_2\xi_2$ and $L_2 = \beta_1\xi_1 +
\beta_2\xi_2$ are independent, then  $\mu_j=\gamma_j*m_{K_j}*E_{g_j}$,
where $\gamma_j \in \Gamma({\mathbb R}^n)$, $K_j$ is a finite subgroup of the group $G$,
$g_j\in G$, $j=1, 2$}.

 {\bf Proof of Theorem 2}.
First we reduce the proof of the theorem to the case
  when $G$ is a torsion group, and then to the case
  when $G$ is a finite group.
 Denote by $Y$ the character group of the group $X$.
 The group  $Y$ is topologically isomorphic to
 the group ${\mathbb R}^n\times H$, where  $H$
 is the character group of the group $G$.
 In order not to complicate the notation we assume that
 $Y={\mathbb R}^n\times H$.
Denote by $x=(t, g)$, where
$t\in {\mathbb R}^n$, $g\in G$, elements of the group $X$, and by
 $y=(s, h)$, where $s\in {\mathbb R}^n$, $h\in H$, elements of the group $Y$.

 1. Consider the factor-group
$X/({\mathbb R}^n\times b_G)$, and denote by  $[(t, g)]$ its elements.
Obviously, $X/({\mathbb R}^n\times b_G)$  is a discrete group.
The character group of the group $X/({\mathbb R}^n\times b_G)$ is topologically
isomorphic to the annihilator $A(Y, {\mathbb R}^n\times b_G)$. It is easy to see that
 $A(Y, {\mathbb R}^n\times b_G)=c_H$. Since  ${\mathbb R}^n$ is a connected component of
 zero  of the group $X$, we have $\alpha({\mathbb R}^n)={\mathbb R}^n$. It follows from
  $b_G=b_X$, that $\alpha(b_G)=b_G$.
 Hence,  $\alpha({\mathbb R}^n\times b_G)={\mathbb R}^n\times b_G$ and
 $\alpha$ induces an automorphism $\hat\alpha$ on the factor-group
$X/({\mathbb R}^n\times b_G)$ by the formula
$\hat\alpha [x]=[\alpha x]$. We note that $\tilde\alpha(c_H)=c_H$ and a topological automorphism
 of the group $c_H$ which is adjoint to $\hat\alpha$, is a restriction of
 $\tilde\alpha$ to $c_H$. Verify that
\begin{equation}\label{1d}
(I+\tilde\alpha_{c_H})(c_H)=c_H.
\end{equation}
Taking into account that for a discrete Abelian group $X$ conditions (\ref{1}) and (\ref{2b})
are equivalent, it suffices to verify that
\begin{equation}\label{2d}
{\rm Ker}(I+\hat\alpha)=\{0\}.
\end{equation}

Take $(t_0, g_0)\in X$, and let $[(t_0, g_0)]\in {\rm Ker}(I+\hat\alpha)$.
Then $[(I+\alpha)(t_0, g_0)]=0$.
Hence, $(I+\alpha)(t_0, g_0)\in {\mathbb R}^n\times b_G$. Since $b_G$ is a torsion
group, we have $k(I+\alpha)(t_0, g_0)\in{\mathbb R}^n$
for some natural $k$. This implies that $(I+\alpha)k(t_0, g_0)\in{\mathbb R}^n$, i.e.
\begin{equation}\label{3d}
(I+\alpha)(kt_0, kg_0)=(t', 0).
\end{equation}
It is obvious that
$(I+\alpha)(\mathbb{R}^n)\subset\mathbb{R}^n$. It follows from $(\ref{1})$
that the restriction of the continuous endomorphism
$I+\alpha$ of the group $X$  to the subgroup
 $\mathbb{R}^n$ is a topological automorphism
 of the group $\mathbb{R}^n$.
Hence,
 \begin{equation}\label{4d}
(t', 0)=(I+\alpha)(\tilde t, 0)
\end{equation}
for an element  $\tilde t\in \mathbb{R}^n$. Taking into account (\ref{1}), it follows from
 (\ref{3d}) and (\ref{4d})
that $kg_0=0$. Hence, $g_0\in b_G$, and it means that
$(t_0, g_0)\in{\mathbb R}^n\times b_G$. This implies that $[(t_0, g_0)]=0$.
So, (\ref{2d}) holds, and thereby (\ref{1d}) is proved.

By Lemma 1,
the symmetry of the conditional distribution of
the linear form
$L_2$ given $L_1$ implies that the characteristic functions
$\hat\mu_j(y)$ satisfy equation
(\ref{2a}). Since $\tilde\alpha(c_H)=c_H$, consider the restriction of
equation (\ref{2a}) to the subgroup $c_H$. Taking into account
(\ref{1d}), apply Lemma 9 to the group $c_H$. We obtain that
 $\hat\mu_j(y)=(x_j,y)$, $y\in c_H$, $j=1, 2$.
By the  extension theorem for characters from a closed subgroup
to the group, we can assume that $x_j\in X$, $j=1, 2$. Substituting
these expressions for
$\hat\mu_j(y)$ to equation (\ref{2a}) and taking into account that $A(X, c_H)={\mathbb R}^n\times b_G$, we get
\begin{equation}\label{12.08.15}
2(x_1+\alpha x_2)\in \mathbb{R}^n\times b_G.
\end{equation}
Since $b_G$ consists of all elements of finite order of
the group $X$, it follows from (\ref{12.08.15}) that
\begin{equation}\label{5d}
x_1+\alpha x_2\in \mathbb{R}^n\times b_G.
\end{equation}
Consider new independent random variables  $\eta_1=\xi_1+\alpha x_2$ and
$\eta_2=\xi_2-x_2$  with values in the group
 $X$. Denote by
  $\lambda_j$  the distribution of the  random variable  $\eta_j$.
  Obviously, ${\lambda_1=\mu_1*E_{\alpha x_2}}$, ${\lambda_2=\mu_2*E_{- x_2}}$.
It is easy to see that the characteristic functions $\hat\lambda_j(y)$
 also satisfy equation (\ref{2a}). Hence, by Lemma 1, the conditional distribution
 of the linear form  $N_2 = \eta_1 +
\alpha\eta_2$ given  $N_1 = \eta_1 +
\eta_2$ is symmetric. Obviously, $\hat\lambda_2(y)=1$, $y\in c_H$. It follows from (\ref{5d})
that  $\hat\lambda_1(y)=1$, $y\in c_H$.
By Lemma 2, this implies that
$\sigma(\lambda_j)\subset A(X, c_H)=\mathbb{R}^n\times b_G$, $j=1, 2$.
Since $\alpha({\mathbb R}^n\times b_G)={\mathbb R}^n\times b_G$,
we reduced the proof of the theorem to the case, when
 $G$ is a torsion group. So, we will assume that.

2. Since $G$ is a torsion group, we have $\alpha(G)=G$. Hence,
$\alpha(t, g)=(\alpha t, \alpha g)$
and $\tilde\alpha(s, h)=(\tilde\alpha s, \tilde\alpha h)$.
Write equation (\ref{2a}) in the form
\begin{equation}\label{1e}
\hat\mu_1(s+s', h+h')\hat\mu_2(s+\tilde\alpha s',h+\tilde\alpha h')$$$$=
\hat\mu_1(s-s', h-h')\hat\mu_2(s-\tilde\alpha s',h-\tilde\alpha h'),
\quad (s, h), (s', h') \in Y.
\end{equation}
Substituting in (\ref{1e}) $s=s'=0$, we get
\begin{equation}\label{2e}
\hat\mu_1(0, h+h')\hat\mu_2(0,h+\tilde\alpha h')=
\hat\mu_1(0, h-h')\hat\mu_2(0,h-\tilde\alpha h'), \quad  h, h' \in H.
\end{equation}
It follows from Theorem 1,   Lemma 1 and
 (\ref{11a}) that solutions of equation (\ref{2e}) are of the form
\begin{equation}\label{3e}
\hat\mu_1(0, h)=(g_1, h)\hat m_K(h), \quad \hat\mu_2(0, h)=(g_2, h)\hat m_K(h),
\quad  h   \in H,
\end{equation}
where $K$ is a finite subgroup of the group $G$, $g_1, g_2\in G$ and $\alpha(K)=K$.
Put $L=A(H, K)$. Substituting (\ref{3e}) in (\ref{2e}), taking into account (\ref{11a}),
and considering the restriction of of the resulting equation to
 $L$, we obtain that ${2(g_1+\alpha g_2)\in K}$.
Since the group $G$ contains no elements of order 2,
 this implies that
\begin{equation}\label{4e}
g_1+\alpha g_2\in K.
\end{equation}
Consider new independent random variables $\eta_1=\xi_1+\alpha g_2$ and
$\eta_2=\xi_2-g_2$  with values in the group
 $X$. Denote by
  $\lambda_j$   the distribution of the  random variable  $\eta_j$. Obviously, ${\lambda_1=\mu_1*E_{\alpha g_2}}$, ${\lambda_2=\mu_2*E_{- g_2}}$.
 It is easy to see that the characteristic functions $\hat\lambda_j(y)$
   also satisfy equation (\ref{2a}). Hence, by Lemma 1, the conditional distribution
 of the linear form $P_2 = \eta_1 +
\alpha\eta_2$ given  $P_1 = \eta_1 +
\eta_2$ is symmetric.
Obviously, $\hat\lambda_2(y)=1$, $y\in L$. It follows from (\ref{5d})
that  $\hat\lambda_1(y)=1$, $y\in L$.
By Lemma 2, this implies that
$\sigma(\lambda_j)\subset A(X, L)=\mathbb{R}^n\times K$, $j=1, 2$.
Since $\alpha({\mathbb R}^n\times K)={\mathbb R}^n\times K$,
we reduced the proof of the theorem to the case, when
 $G$ is a finite group, and we will assume that.

3. By Lemma 3, the linear forms
$M_1=(I+\alpha)\xi_1+2\alpha\xi_2$ и
$M_2=2\xi_1+(I+\alpha)\xi_2$ are independent. Since $G$ is a finite group
containing no elements of order 2, it follows from
(\ref{1}) that the coefficients of the
linear forms $M_1$ and $M_2$ are topological automorphisms
of the group $X$.
By Lemma  11, this implies that $\mu_j=\gamma_j*m_{K_j}*E_{g_j}$,
where $\gamma_j \in \Gamma({\mathbb R}^n)$, $K_j$ is a finite subgroup of
  $G$,
$g_j\in G$, $j=1, 2$. The equality $K_1=K_2=K$ follows from (\ref{3e}).  Theorem 2 is proved.

We note that Theorem 2 was proved in  \cite{Fe7} for the groups $X$ of the form
$X={\mathbb R}^n\times G$, where $n \ge 0,$ and $G$ is
a discrete Abelian group containing no elements
of order $2$ and such that its torsion part is a finite group.
We also note that if a topological automorphism $\alpha$ satisfies conditions (\ref{1a}), then   the statement of Theorem 2 was proved
 in \cite[Theorem 2 and Remark 4]{Fe3}, see also
 \cite[Remark 17.27]{Fe5a}. Reasoning as in Remark  4, we see that Theorem 2 fails if we omit condition (\ref{1}).

We prove now the following generalization of Theorem  1.

{\bf Theorem 3.}
 {\it  Let $X$ be a discrete Abelian group, , $F$ be the  $2$-component of the group $X$, $G$ be the  subgroup generated by all elements of add order of the group $X$.
  Let  $\alpha$ be an automorphism of
 $X$ satisfying condition $(\ref{1})$.
Let
  $\xi_1$ and  $\xi_2$ be independent random variables with values in the group
       $X$  and distributions $\mu_1$ and $\mu_2$.
If the conditional distribution of the linear form
$L_2 = \xi_1 + \alpha\xi_2$ given  $L_1 = \xi_1 +
\xi_2$  is symmetric, then
$\mu_j=\rho_j*m_K*E_{x_j}$, where $\sigma(\rho_j)\subset F$,  $K$ is a finite subgroup of
 $G$, $x_j\in X$, $j=1, 2$.}

 To prove Theorem 3 we need the following lemma.

 {\bf Lemma 12}.  {\it Let  $X$ be a  discrete torsion Abelian group,
 $\alpha$ be an automorphism of
 the group $X$ satisfying condition $(\ref{1})$.
Let
  $\xi_1$ and  $\xi_2$ be independent random variables with values in the group
       $X$  and distributions $\mu_1$ and $\mu_2$.
If the conditional distribution of the linear form
$L_2 = \xi_1 + \alpha\xi_2$ given  $L_1 = \xi_1 +
\xi_2$  is symmetric, then
$\mu_j$ are supported in  a subgroup generated by  $X_{(2)}$
and a finite subgroup, and hence, in a subgroup $X_{(n)}$ for some  $n$}.

{\bf Proof}.  First we prove the lemma supposing
 that $\hat \mu_j(y)\ge 0$, $j=1, 2$. Reasoning as in the proof of Lemma 5
 and retaining notation of Lemma 5  we obtain that
 \begin{equation}\label{13.08.15}
\hat\mu_1(y)=\hat\mu_1(y)=1, \quad y\in B,
\end{equation}
but in contrast to Lemma 5,
if the group $X$ contains elements of order 2, then
the subgroup $B$ does not need to be open.

We have $W=V^{(2)}\cap(I+\tilde\alpha)(V)\cap (\tilde\alpha(V))^{(2)}$.
Put $L=V\cap(I+\tilde\alpha)(V)\cap\tilde\alpha(V)$, $M=(I+\tilde\alpha)(L)$.
By (\ref{2b}), the continuous endomorphism $I+\tilde\alpha$
is open. This implies that $L$, and hence $M$  are open subgroups in
$Y$. It is obvious that $W\supset L^{(2)}$. Hence,
$B\supset M^{(2)}$. Put $G=A(X, M^{(2)})$. By Lemma 2, it follows from
(\ref{13.08.15}) that $\sigma(\mu_j)\subset A(X, B
)\subset G$, $j=1, 2$.

We have $G=\{x\in X: (x, 2y)=1$ for all
$y\in M\}=\{x\in X: (2x, y)=1$ for all
$y\in M\}=\{x\in X: 2x\in A(X, M)\}=f_2^{-1}(A(X, M)).$
Note that ${\rm Ker}f_2=X_{(2)}$.
Since $M$  is an open subgroup in
 $Y$, the annihilator $A(X, M)$ is a compact subgroup in $X$, and hence
 a finite subgroup in $X$. This implies that the group $G=f_2^{-1}(A(X, M))$
is generated by  $X_{(2)}$
and a finite subgroup of $X$.

We proved the lemma, assuming that $\hat \mu_j(y)\ge 0$, $j=1, 2$.
Get rid of this restriction.
Put
$\nu_j = \mu_j
* \bar \mu_j$. Then
 $\hat \nu_j(y) = |\hat \mu_j(y)|^2 \ge 0,$   $y \in Y$.
Obviously, the characteristic functions $\hat \nu_j(y)$
also satisfy equation (\ref{2a}). Denote by $\tilde\xi_j$ the independent
random variables with values in the group $X$ and distributions $\nu_j$,
$j=1, 2$. By Lemma 1, the conditional distribution of the linear form
 $\tilde L_2 = \tilde\xi_1 + \alpha\tilde\xi_2$ given $\tilde L_1 = \tilde\xi_1 + \tilde\xi_2$ is symmetric.
As has been proved above
$ \nu_j $ are supported in   a subgroup $S$, generated by
 $X_{(2)}$ and a finite subgroup.
It follows from this that each distribution  $\mu_j$
is supported in a set $x_j+S$ for some
 $x_j\in X$. Since $X$ is a torsion group, each of the elements
$x_j$ has a finite order. Hence, the subgroup generated by
 $S$ and $x_j$ also has the required form.
Lemma 12 is proved.

{\bf Proof of Theorem 3}.  Since $\alpha(b_X)=b_X$, we can apply Lemma 10 and
 assume that $X$ is a torsion group. Since $\alpha(X_{(n)})=X_{(n)}$, by Lemma 12, we can assume that  $X=X_{(n)}$ for some
$n$. Since  $X$  is a torsion group,   $X$
is
a weak direct product of its $p$-torsion subgroups
(\cite[Theorem 8.4]{Fu1}). This implies that, $X=F\times G$.
Denote by $Y$ the character group of the group $X$.
 The group  $Y$ is topologically isomorphic to the group $L\times H$,
 where $L$ is the character group of $F$, and $H$ is
 the character group of $G$.
In order not to complicate the notation we assume that $Y=L\times H$.
Since $\alpha(F)=F$ and $\alpha(G)=G$, we have $\tilde\alpha(L)=L$, $\tilde\alpha(H)=H$.
 Denote by $y=(l, h)$,
where $l\in L$, $h\in H$, elements of the group $Y$.

Reasoning as in the proof of item 2 of Theorem 2, and considering
the group $F$ instead of $\mathbb{R}^n$, we reduce the proof of the theorem to
the case, when $G$  is a finite group, and
\begin{equation}\label{21.04.15.3}
\hat\mu_1(0, h)=\hat\mu_2(0, h)=\hat m_G(h),
\quad  h   \in H.
\end{equation}
By Lemma 1, the symmetry of the conditional distribution of the linear form
$L_2$ given $L_1$ implies that the characteristic functions
 $\hat\mu_j(y)$ satisfy equation
(\ref{2a}).
Substitute in (\ref{2a}) $u=v=(l, h)$. We get
\begin{equation}\label{21.04.15.2}
\hat\mu_1(2(l, h))\hat\mu_2((I+\tilde\alpha)(l, h))=
\hat\mu_2((I-\tilde\alpha)(l, h)), \quad (l, h) \in Y.
\end{equation}
Substitute in (\ref{21.04.15.2}) $l=0$, $h\ne 0$. Since $2h\ne 0$, it follows from
 (\ref{21.04.15.3}) that the left-hand side of the resulting equation
 is equal to zero. Thus    $\hat\mu_2(0, (I-\tilde\alpha_H)h)=0$.
Hence,  if  $h\ne 0$, then $(I-\tilde\alpha_H)h\ne 0$, i.e.
${\rm Ker}(I-\tilde\alpha_H)=\{0\}$.
Since $H$ is a finite group, this implies that
\begin{equation}\label{21.04.15.6}
I-\tilde\alpha_H\in{\rm Aut}(H).
\end{equation}

Obviously, the automorphism  $\alpha_F$  satisfies the condition
\begin{equation}\label{21.04.15.7}
{\rm Ker}(I+\alpha_F)=\{0\}.
\end{equation}
It is easy to see that (\ref{21.04.15.7})   is equivalent to the condition
\begin{equation}\label{21.04.15.4}
{\rm Ker}(I-\alpha_F)=\{0\}.
\end{equation}
Indeed, assume that (\ref{21.04.15.7}) is fulfilled and
$x\in{\rm Ker}(I-\alpha_F),$ $x\ne 0$.
We can assume without loss of generality that $2x=0$. Then the equality $(I+\alpha_F)x=(I-\alpha_F)x+2\alpha_F x=0$ implies the contradiction. Reasoning similarly, we get that (\ref{21.04.15.4})
implies (\ref{21.04.15.7}).

Consider the factor group
$X/X_{(2)}$.
Its character group is topologically isomorphic to
the annihilator $A(Y, X_{(2)})$. Obviously,
 $A(Y, X_{(2)})=Y^{(2)}$.
 It follows from $\alpha(X_{(2)})=X_{(2)}$ that
 $\alpha$ induces an automorphism $\hat\alpha$ of the
 factor-group $X/X_{(2)}$ by the formula $\hat\alpha [x]=[\alpha x]$.
 In so doing, a topological automorphism
 of the group $Y^{(2)}$  which is adjoint to $\hat\alpha$
is the restriction of  $\tilde\alpha$ to $Y^{(2)}$.

Obviously, $\tilde\alpha(Y^{(2)})=Y^{(2)}$.
Verify that (\ref{3b}) is fulfilled.
Take $x_0\in X$ such that $[x_0]\in {\rm Ker}(I+\hat\alpha)$. Then $[(I+\alpha)x_0]=0$, and hence
$(I+\alpha)x_0\in X_{(2)}$, i.e. $2(I+\alpha)x_0=0$. It follows from
 (\ref{1}) that $2x_0=0$, i.e.   $[x_0]=0$.
 Thus  (\ref{3b})  holds true.

Since we assume that
$X=X_{(n)}$, the group $X$ is bounded. Hence, the group $F$ is also bounded.
Denote by $k_X$ the least nonnegative integer such
that
$F_{(2^{k_X})}=\{0\}$. If $k_X=0$, then $X=G$, and the statement of
the theorem follows from Theorem 1. We prove by induction
that if the theorem holds true for the groups $X$ satisfying
the condition $k_X=m-1$, then it is valid
for the groups $X$ satisfying
the condition  $k_X=m$.

So, let $X$ be a group such that
$k_X=m$. Put  $\hat\mu_1(l, 0)=a_1(l)$,   $\hat\mu_2(l, 0)=a_2(l)$ and prove that
\begin{equation}\label{21.04.15.1}
\hat\mu_1(l, h)= \begin{cases}a_1(l), & \text{\ if\ }h=0,\\ 0, &   \text{\ if\ }h\ne 0,
\\
\end{cases} \quad\quad\quad
\hat\mu_2(l, h)= \begin{cases}a_2(l), & \text{\ if\ }h=0,\\ 0, &   \text{\ if\ }h\ne 0.
\\
\end{cases}
\end{equation}
The statement of the theorem easily follows from this.

The restrictions of the characteristic functions $\hat\mu_j(l, h)$ to the subgroup $Y^{(2)}$
are the characteristic functions of some independent random variables
 $\eta_1$  and $\eta_2$ with values in the factor-group $X/X_{(2)}$. By Lemma 1,
the conditional distribution of the linear form  $M_2 = \eta_1 + \hat\alpha\eta_2$ given $M_1 = \eta_1 + \eta_2$ is symmetric.
  As has been noted above, (\ref{3b}) is fulfilled. Obviously, $k_{X/X_{(2)}}=m-1$.
Then by induction hypothesis  (\ref{21.04.15.1}) holds true for $y\in Y^{(2)}$.
Let $(l, h) \in Y$, and $h\ne 0$. Since $2(l, h)\in Y^{(2)}$ and $2h\ne 0$, the left-hand side  of equation (\ref{21.04.15.2}) is equal to zero.
 Hence, $\hat\mu_2((I-\tilde\alpha)(l, h))=0$
for all $(l, h) \in Y$, $h\ne 0$.
We have, $(I-\tilde\alpha)(l, h)=((I-\tilde\alpha_L) l, (I-\tilde\alpha_H )h).$
It follows from (\ref{21.04.15.4}) that $(I-\tilde\alpha_L)(L)=L$.
Taking into account (\ref{21.04.15.6}), we get that representation (\ref{21.04.15.1}) for the function $\hat\mu_2(l, h)$  holds true for all $(l, h)\in Y$.
Substituting in (\ref{2a}) $u=(l, h)$, $v=\tilde\alpha^{-1}(l, h)$ and reasoning similarly, we
obtain the required representation for the function
$\hat\mu_1(l, h)$.  Theorem 3 is proved.

Note that if an automorphism
 $\alpha$ satisfies conditions
 (\ref{1a}), then  Theorem 3 is proved
 in  \cite{My2}, see also \cite[Theorem 17.33]{Fe5}.

Condition (\ref{1}) is necessary and sufficient for Theorem 1 and 2. As to Theorem 3
the situation is more complicated.  The reasoning  similar to that used in Remark 4
shows that the condition
 \begin{equation}\label{14.11.1}
{\rm Ker}(I+\alpha)\subset X_2
\end{equation}
is necessary if we want that Theorem 3 be valid.
Below we consider the case, when condition (\ref{14.11.1}) holds, but condition (\ref{1}) fails, i.e.
 \begin{equation}\label{14.11.2}
{\rm Ker}(I+\alpha)\ne\{0\}.
\end{equation}
We complement Theorem 3 by the following statement, which shows that,
 generally speaking, Theorem 3 fails if we change condition  (\ref{1}) for   (\ref{14.11.1}).

{\bf Proposition  2}. {\it Let $X$ be a discrete Abelian group, $F$ be the $2$-torsion subgroup of $X$,   $G$ be the subgroup generated by all elements of add order of  $X$.
 Let  $\alpha$ be an automorphism of the group $X$ satisfying condition   $(\ref{14.11.1})$.
 Then the following statements hold.

  $1$. Assume that the only finite  subgroup $K$ of the group
   $G$, satisfying the condition  $(I-\alpha)(K)=K$, is $K=\{0\}$.
   Let
  $\xi_1$ and  $\xi_2$ be independent random variables with values in the group
       $X$  and distributions $\mu_1$ and $\mu_2$.
If the conditional distribution of the linear form
$L_2 = \xi_1 + \alpha\xi_2$ given  $L_1 = \xi_1 +
\xi_2$  is symmetric, then
    $\mu_j=\rho_j*E_{x_j}$, where $\sigma(\rho_j)\subset F$,  $x_j\in X$, $j=1, 2$.

$2$. Assume that there exists a non-zero finite subgroup $K_0$ of the group $G$, satisfying the condition $(I-\alpha)(K_0)=K_0$. Assume that the automorphism $\alpha$ satisfies the condition
  $(\ref{14.11.2})$.
 Then there exist independent  identically distributed random variables
$\xi_1$ and $\xi_2$ with values in the group $X$ and distribution $\mu$ such that
$\mu$ can not be represented in the form  $\mu=\rho*m_K*E_{x}$, where
 $\sigma(\rho)\subset F$,  $K$ is a finite subgroup of the group
 $G$, $x\in X$, whereas the conditional distribution of the linear form
 $L_2 = \xi_1 + \alpha\xi_2$ given $L_1 = \xi_1 + \xi_2$ is symmetric.}

{\bf Proof}. 1. By Lemma 10, we can assume that $X$ is a torsion group.
Then we have $X=F\times G$.
Denote by $Y$ the character group of the group $X$.
 The group  $Y$ is topologically isomorphic to the group $L\times H$,
 where $L$ is the character group of $F$, and $H$ is
 the character group of $G$.
In order not to complicate the notation we assume that $Y=L\times H$.
Since $\alpha(F)=F$ and $\alpha(G)=G$, we have $\tilde\alpha(L)=L$, $\tilde\alpha(H)=H$.
 Denote by $y=(l, h)$,
where $l\in L$, $h\in H$, elements of the group $Y$.
By Lemma 1, the characteristic functions $\hat\mu_j(y)$ satisfy equation  (\ref{2a}). Put
$\nu_j = \mu_j
* \bar \mu_j$. Then
 $\hat \nu_j(y) = |\hat \mu_j(y)|^2 \ge 0,$   $y \in Y$. Obviously,
 the characteristic functions  $\hat \nu_j(y)$
also satisfy equation (\ref{2a}), which takes the form
\begin{equation}\label{14.11.8}
\hat\nu_1(l+l', h+h')\hat\nu_2(l+\tilde\alpha_L l',h+\tilde\alpha_H h')$$$$=
\hat\nu_1(l-l', h-h')\hat\nu_2(l-\tilde\alpha_L l',h-\tilde\alpha_H h'),
\quad  l, l'
\in L, \ h, h'
\in H.
\end{equation}
Putting here $l=l'=0$,    we obtain
\begin{equation}\label{18.11.1}
\hat\nu_1(0, h+h')\hat\nu_2(0,h+\tilde\alpha_H h')=
\hat\nu_1(0, h-h')\hat\nu_2(0,h-\tilde\alpha_H h'),
 \quad h, h'
\in H.
\end{equation}
  It follows from  (\ref{14.11.1}) that ${\rm Ker}(I+\alpha_G)=\{0\}$. This implies by
 Lemma 1 and Theorem 1, applying to the group  $G$  that
 $\hat\nu_1(0, h)=\hat\nu_2(0, h)=\hat m_K(0, h)$, $h\in H$, where
 $K$ is a finite subgroup of $G$. We have by Lemma 1 and Proposition 1 that $(I-\alpha)(K)=K$.
 It follows from the condition of Proposition 2 that  $K=\{0\}$. Hence, $\hat\nu_1(0, h)=\hat\nu_2(0, h)=1$, $h\in H$. By Lemma 2, this implies that
  $\sigma(\nu_j)\subset F$, and hence, $\mu_j=\rho_j*E_{x_j}$, where $\sigma(\rho_j)\subset F$,  $x_j\in X$, $j=1, 2$. Thus, we proved statement 1.

2. It is easy to see that $\alpha(K_0)=K_0$. Consider the subgroup  $T=F_{(2)}\times K_0$ of the group $X$. Obviously, $\alpha(T)=T$, i.e. the restriction of
$\alpha$ to $T$ is an automorphism of the group $T$. Denote by $\alpha_T$ this restriction.
Then the group $T$ and the automorphism $\alpha_T$ satisfy the conditions of statement 2.
Taking this into account we can prove statement 2, assuming that $X=F\times G$,
\begin{equation}\label{14.11.4}
F=F_{(2)},
\end{equation}
and $G$ is a finite group such that $(I-\alpha_G)(G)=G$.
Since $G$ is a finite group, it means that
\begin{equation}\label{17.11.1}
I-\alpha_G\in {\rm Aut}(G).
\end{equation}
It follows from (\ref{14.11.1}) and (\ref{14.11.2}) that
 ${\rm Ker}(I+\alpha_F)\ne\{0\}$, and hence, $(I+\tilde\alpha_L)(L)=A(L, {\rm Ker}(I+\alpha_F))\ne L$.
 Take $l_0\notin (I+\tilde\alpha_L)(L)$. Since $F$ is a torsion group,
  $L$ is a compact and totally disconnected group.
 For this reason there exists an open in $L$  subgroup $U$
 such that $U\cap (l_0+U)=\emptyset$ and $(l_0+U)\cap (I+\tilde\alpha_L)(L)=\emptyset$.
  Consider the factor-group $A=Y/U\cong L/U\times H$,
   and denote by $[y]=([l], h)$, where $[l]\in L/U$, $h\in H$,
 its elements. Denote by $n$
  the number of elements of the  group $H$. Consider on the group
$A$ the function
$$
f([y])=f([l], h)= \begin{cases}1, & \text{\ if\ }[l]=0, \ h=0,\\ 0, &   \text{\ if\ }[l]=0, \ h\ne 0,
\\
{1\over n}, &   \text{\ if\ }[l]=[l_0], \ h\in H,\\
0, &   \text{\ if\ }[l]\notin\{0, [l_0]\}, \ h\in H.
\end{cases}
$$
Put $B=A(F, U)\times G$.
  By the duality theorem the groups
 $A$ and $B$ are the character groups of each other.
 Put
$$p(x)=1+\sum\limits_{[y]\in A, \ [y]\ne 0}f([y])\overline{(x, [y])}, \quad x\in B.$$
It is obvious that $p(x)\ge 0$ and $\int_Bp(x)dm_B(x)=1$. Let $\mu$ be the distribution
on the group $B$ with density $p(x)$ with resect to $m_B$. We have
$\hat\mu([y])=f([y])$, $[y]\in A$. If we consider $\mu$ as a distribution on the group $X$, then the characteristic function $\hat\mu(y)$ is of the form
\begin{equation}\label{14.11.5}
\hat\mu(y)=\hat\mu(l, h)= \begin{cases}1, & \text{\ if\ }l\in U, \ h=0,\\ 0, &   \text{\ if\ }l\in U, \ h\ne 0,
\\
{1\over n}, &   \text{\ if\ }l\in l_0+U, \ h\in H,\\
0, &   \text{\ if\ }l\notin U\cup (l_0+U), \ h\in H.
\end{cases}
\end{equation}
 Obviously,  $\mu$ can not be represented in the form
 $\mu=\rho*m_K*E_{x}$, where $\sigma(\rho)\subset F$,  $K$ is a finite subgroup of the group
 $G$, $x\in X$. Let $\xi_1$ and $\xi_2$  be independent  identically distributed random variables
 with values in the group $X$ and distribution $\mu$.  It follows from
 (\ref{14.11.4}) that $L_{(2)}=L$. Hence $l=-l$, $l\in L$.  We shall verify that
 the characteristic function $f(l, h)=\hat\mu(l, h)$ satisfies the equation
\begin{equation}\label{14.11.6}
f(l+l', h+h') f(l+\tilde\alpha_L l', h+\tilde\alpha_H h')= f(l+l', h-h') f(l+\tilde\alpha_L l', h-\tilde\alpha_H h'), \ l, l'
\in L, \ h, h'
\in H.
\end{equation}
Then, by Lemma 1,  the conditional distribution of the linear form
 $L_2 = \xi_1 + \alpha\xi_2$ given $L_1 = \xi_1 + \xi_2$ is symmetric.
The following cases exhaust all possibilities for $l$ and $l'$.

$(i)$. $l+l'\in U$ and $l+\tilde\alpha_L l'\in U$. It follows from
(\ref{14.11.5}) that  the left-hand side and the right-hand side of equation (\ref{14.11.6}) can take
values either 1 or 0. Assume that the left-hand side
of equation (\ref{14.11.6}) is equal to 1. This implies that
 $h+h'=0$ and $h+\tilde\alpha_H h'=0$. Hence,
\begin{equation}\label{17.11.2}
(I-\tilde\alpha_H) h'=0.
\end{equation}
It follows from (\ref{17.11.1}) that $I-\tilde\alpha_H \in {\rm Aut}(H)$, and  (\ref{17.11.2})
implies that $h'=0$, and hence $h=0$. It means that the right-hand side
of equation (\ref{14.11.6}) is also equal to 1. Reasoning similarly, we verify that
if right-hand side
of equation (\ref{14.11.6}) is equal to 1, then the left-hand side
of equation (\ref{14.11.6}) is equal to 1. Thus, in this case both sides of equation (\ref{14.11.6}) are equal.

 $(ii)$. $l+l'\in l_0+ U$ and
$l+\tilde\alpha_L l'\in l_0+U$. It follows from (\ref{14.11.5}) that both sides of equation (\ref{14.11.6}) are ${1\over n^2}={1\over n^2}$.

$(iii)$.  $l+l'\in   U$ and  $l+\tilde\alpha_L l'\in l_0+ U$. This implies that
$(I+\tilde\alpha_L)l'\in l_0+ U$. But this contradicts the fact that
$(l_0+U)\cap (I+\tilde\alpha_L)(L)=\emptyset$. Thus, this case is impossible.

$(iv)$. $l+l'\in l_0+ U$ and  $l+\tilde\alpha_L l'\in   U$. The reasoning is similar to the case $(iii)$.

$(v)$. $l+l'\notin l_0+ U$ and  $l+l'\notin U$. Then in view of  (\ref{14.11.5}) we have $f(l+l', h+h')= f(l+l', h-h')=0$, and both sides of equation (\ref{14.11.6}) are zero.

$(vi)$. $l+\tilde\alpha_L l'\notin l_0+ U$ and  $l+\tilde\alpha_L l'\notin U$. Then in view of (\ref{14.11.5}) we have $f(l+\tilde\alpha_L l', h+\tilde\alpha_H h')= f(l+\tilde\alpha_L l', h- \tilde\alpha_H h')=0$, and both sides of equation (\ref{14.11.6}) are zero.

 Proposition 2 is proved.

\vskip 2 cm

\noindent B. Verkin Institute for Low Temperature
Physics and Engineering of the
National Academy of Sciences of Ukraine,
Kharkiv, Ukraine; feldman@ilt.kharkov.ua

\end{document}